\documentclass[11pt,letterpaper,reqno]{amsart}

\usepackage{enumitem}
\usepackage{amsfonts}
\usepackage{amssymb}
\usepackage{amsmath}
\usepackage{amsthm}
\usepackage{mathtools}
\usepackage{amstext}
\usepackage{tikz-cd}
\usepackage[left=1.3in,right=1.3in,top=1.2in, bottom=1.3in]{geometry}

\usepackage{graphicx}
\definecolor{cit}{HTML}{117733}
\definecolor{lin}{HTML}{0003d6}
\definecolor{gree}{HTML}{6dc066}
\usepackage{hyperref}
\hypersetup{bookmarksnumbered,colorlinks=true,linkcolor=lin,citecolor=cit,urlcolor=black}

\usepackage[section]{placeins}
\usepackage{caption}

\DeclareMathOperator{\II}{II}
\DeclareMathOperator{\Int}{int}

\DeclareMathOperator{\PGL}{PGL}
\DeclareMathOperator{\kum}{Kum}

\DeclareMathOperator{\Aut}{Aut}

\DeclareMathOperator{\bl}{Bl}

\DeclareMathOperator{\amp}{Amp}
\DeclareMathOperator{\Nef}{Nef}

\DeclareMathOperator{\Pic}{Pic}
\DeclareMathOperator{\NE}{NE}
\DeclareMathOperator{\NEb}{\overline{\NE}}

\DeclareMathOperator{\id}{id}

\newcommand{\Z}{\mathbb{Z}}

\newcommand{\R}{\mathbb{R}}
\newcommand{\PP}{\mathbb{P}}

\newcommand{\Q}{\mathbb{Q}}

\newcommand{\vv}[1]{\lvert#1\rvert}
\newcommand{\vr}[1]{\langle#1\rangle}
\newcommand{\ra}{\to}

\newcommand{\mt}{\mapsto}

\newcommand{\ik}{^{-1}}

\newcommand{\txt}[1]{\text{ #1}} 
\theoremstyle{plain}
\newtheorem{theorem}{Theorem}[section]

\newtheorem{prop}[theorem]{Proposition}
\newtheorem{lemma}[theorem]{Lemma}
\newtheorem{conj}[theorem]{Conjecture}

\theoremstyle{definition}
\newtheorem{defi}[theorem]{Definition}

\newtheorem{remark}[theorem]{Remark}

\newcommand{\pf}{{\em Proof}}
\newcommand{\pfof}[1]{{\em Proof of #1}}

\begin{document}
\title[Pseudo-automorphisms and Kummer surfaces]{Pseudo-automorphisms of rational threefolds and Kummer surfaces}
\author{Zhuang He}
\email{he.zhu@northeastern.edu}
\date{\today}

\begin{abstract}
	Kummer surfaces are special quartic surfaces that admit $16$ nodes. The automorphisms of K3 Kummer surfaces are rich and complicated. Based on the results of Keum and Kond\=o, and as a continuation of the recent result by He and Yang, we lift $45$ classically known automorphisms of Kummer surfaces to pseudo-automorphisms of a threefold, the blow-up of $\mathbb{P}^3$ along $6$ points and $15$ lines. We give a description of a fundamental domain of the group generated by all the known pseudo-automorphisms on this threefold.
\end{abstract}
\maketitle
\section{Introduction}
The theory of automorphisms of higher dimensional algebraic projective varieties has a unique favor from the curves and surfaces.
In dimensions at least $3$, automorphisms are much rarer than birational automorphisms. Some special birational automorphisms are of particular interest, that are called {\em pseudo-automorphisms}. These are birational automorphisms that is an isomorphism outside of a locus of codimension at least 2. Equivalently, by \cite{Cantat2017}, pseudo-automorphisms are birational maps that neither it and its inverse contracts any divisors.

On rational varieties, or for a specific example, blow-ups of projective spaces, it is a non-trivial question to describe new pseudo-automorphisms. One interesting case is given by \cite{HY2019}, on the blow-up of $\mathbb{P}^3$ of six very general points and the $15$ lines through them. The anticanonical divisor of this space, which we denote by $X$, turns out to be a Jacobian Kummer surface. 

Kummer surfaces are quartic surfaces in $\mathbb{P}^3$ that admit $16$ nodes, which is the maximal possible number for a quartic. Resolving the $16$ nodes gives a smooth K3 surface $S$, of Picard number at least $17$. Although the geometry on Kummer surfaces has been classical studied for over a century (\cite{Hutchinson1901, Hudson1905, Hudson1990}), the automorphism group of a general Jacobian Kummer surface is not completely known until the 1990s. Besides the ``classical" automorphisms (full description in Section \ref{Cinv}), which are all involutions themselves, Keum \cite{Keum1997} first constructed an automorphism on Jacobian Kummer surfaces of Picard number $17$, which is of infinite order. Kond\=o later \cite{Kondo1998} described Keum's automorphisms as a family of $192$ automorphisms, interpreted via the theory of Leech roots. Then Kond\=o proved that Keum's 192 automorphisms and the classical automorphisms together generate the automorphism group of Kummer surfaces.

New realizations of Keum's automorphism emerge from from the observation that Jacobian Kummer surfaces are the unique anticanonical section of the threefold $X$.
In \cite{HY2019}, the authors showed that for a fixed pair $(X,S)$, $120$ out of the $192$ Keum's automorphisms lift to pseudo-automorphisms of $X$. The discovery of those pseudo-automorphisms, which are of infinite order, implies that the effective cones of $X$ are not rational polyhedral, and in particular, proved that $X$ is not a Mori Dream Space.

In this paper, we continue the approach in \cite{HY2019}, and find more pseudo-automorphisms of $X$. We also provide some insights to the full description of the pseudo-automorphism group $PsAut(X)$ of $X$.

First of all, every pseudo-automorphism of $X$ restricts to an birational automorphism of $S$. Since $S$ is a K3 surface, every birational automorphism of $S$ is in fact regular. Hence there is a group homomorphism $u: PsAut(X)\ra \Aut(S)$. To begin with, we proved a surprising result:
\begin{prop}(Proposition \ref{injhom})
The group homomorphism $u:PsAut(X)\ra \Aut(S)$ is injective.
\end{prop}
In other words, every pseudo-automorphism of $X$ is in fact a Kummer automorphism. It is then an interesting question to decide $PsAut(X)$ as a subgroup of $\Aut(S$). A necessary condition for a Kummer automorphism to be in $PsAut(X)$, is that they must fix the class $K_X\mid_S$. However, we constructed some non-examples via computer-assisted calculation, that lie in the stabilizer of $K_X\mid_S$, but do not lift to $X$. 

To seek more pseudo-automorphisms, we imposed an extra condition, that the pseudo-automorphism also fixes the unique rational normal curve $R$ through the six points blown-up. Intuitively this makes sense, as the choice of the $R$ among the $16$ otherwise symmetric tropes on a Kummer surface, decides the embedding of $S$ as the anticanonical locus of $X$. Besides, the $120$ Keum's pseudo-automorphisms in \cite{HY2019} are exactly those among the original $192$ ones that fix the curve $R$. Under this assumption, we discovered another family of $45$ involutions of $S$, that lift to pseudo-automorphisms of $X$:

\begin{theorem}
There exist $45$ pseudo-automorphisms of $X$, all are involutions, whose restrictions to the Jacobian Kummer surface $S$ are the composition of a Hutchinson-G\"opel involution and a translation automorphism. Each involution is indexed by a quadruple of nodes called G\"opel Tetrads.
\label{main01}
\end{theorem}

We call the $45$ pseudo-automorphisms above the {\em Hutchinson-G\"opel-type} (HG-type for short) involutions of $X$. It turns out that these $45$ HG-type pseudo-automorphisms, together with the $120$ Keum's pseudo-automorphisms, already give us all the pseudo-automorphisms fixing $R$:
\begin{theorem}\label{main02}
	Let $G$ be the subgroup of $\Aut(S)$ generated by the $45$ HG-type and $120$ Keum's automorphisms.
	Let $B$ be the rank $2$ lattice in the N\'{e}ron-Severi lattice $NS(S)$ spanned by $b:=K_X\mid_S$ and the rational normal curve $R$. Let $A=B^\perp$ be the rank $15$ orthogonal complement of $B$ in $NS(S)$. Then we have:
	\begin{enumerate}
		\item the stabilizer subgroup $\Aut(S)_A$ equals $G$;
		\item In particular,
			\[G=\Aut(S)_A=\Aut(S)_b\cap \Aut(S)_R.\]
	\end{enumerate}
\end{theorem}

Finally we shed light upon the question of fully deciding the pseudo-automorphism group.

\begin{theorem}
	\label{main03}
	\label{Gstab}
	There is a fundamental domain $\Omega$ of the action of $G$ on  $\Nef(S)\cap A$. This $\Omega$ is rational polyhedral. Among the faces of $\Omega$, $45+120$ faces one-to-one correspond to the HG-type and Keum's automophisms in $G$. The remaining faces are cut out by $(-2)$ curve classes.
\end{theorem}

\begin{conj}
	For every $g\in PsAut(X)$, $g(R)=R$. In particular, $PsAut(X)=G$.
\end{conj}

Structure of this paper: Section 2 reviews the geometry and automorphisms of Jacobian Kummer surfaces. Section 3 constructs the HG-type pseudo-automorphisms on $X$, and proves Theorem \ref{main01}. Section 4 reviews Kond\={o}'s construction of a fundamental domain of $\Aut(S)$ acting on the positive cone of $S$, and then proves Theorems \ref{main02} and \ref{main03}. 

{\bfseries Acknowledgement} The author thanks Ana-Maria Castravet and Antonio Laface for talks and suggestions. The main idea of this paper was first formed during the Oberwolfach workshop ``Toric geometry" (2019). The author thanks the organizers for the event.

\section{Kummer surfaces}
\subsection{Nodes and tropes}
We review briefly the geometry of a Jacobian Kummer surface $S$ and its Picard lattice. Let $C$ be a smooth genus $2$ curve.
The singular Kummer surface \[\kum(J(C)):=J(C)/\iota\] is the quotient of the Jacobian of $C$ by the involution $\iota:a\mt -a$. $\kum(J(C))$ embeds in $\PP^3$ as a quartic with sixteen nodes $N_\alpha$, where $\alpha$ runs through the sixteen $2$-torsion points on $J(C)$. There are $16$ hyperplanes in $\PP^3$ tangent to $\kum(J(C))$ along a smooth rational curve $T_\beta$, each through six nodes in $\PP^3$. Those sixteen $T_\beta$ are called the tropes on $\kum(J(C))$. Define $S$ to be the minimal desingularization of $\kum(J(C))$, so each $N_\alpha$ gives an exceptional divisor, which are smooth rational curves, which we also write as $N_\alpha$. $S$ is called the Jacobian Kummer surface.
The incidental relations between the tropes and nodes of $S$ form a combinatorial object, called the $(16,6)$-configuration, where each node lies exactly on six tropes, and each trope passes through exactly six nodes. 

We follow one of the notations in \cite{Dolgachev2002}. Here $p_1,\cdots,p_6$ are the six Weierstrass points on $C$. The map $C\ra J(C)$, $p\mt [p-p_6]$ sends $\{p_i\}$ to the Theta divisor $\Theta_0$, and $J(C)_2$ consists of $0$ and $\{ij$, $1\leq j,k\leq 6\}$. The addition on $J(C)_2$ is given by $0$ being identity, $ij+jk=ik$, and $ij+kl=mn$, for $\{i,j,k,l,m,n\}=\{1,2,3,4,5,6\}$. Then both the nodes and tropes can be indexed by $J(C)_2$. We define $I(T_\beta)$ to be the set of six nodes $T_\beta$ passes through. Then $I(T_0)=\Theta_0$, and $I(T_{ij})$ is exactly the Theta divisor $\Theta_{ij}$, which equals $\Theta_0$ translated by the node $N_{ij}$.
Assume \[\{i,j,k,l,m\}=\{1,2,3,4,5\}.\] Then
\begin{align}
	I(T_0)&=\{0=66,16,26,36,46,56\},\\
	I(T_{i6})&=\{0,i6,ij,ik,il,im\},\\
	I(T_{ij})&=\{i6,j6,ij,lm,ln,mn\}.
\end{align}
Now the Picard lattice $\Pic(S)$ coincides with the N\'{e}ron-Severi lattice $NS(S)$. For $C$ very general, $NS(S)$ has rank $17$, and $NS(S)$ is spanned by $\{N_\alpha, T_\beta\}$ over $\Z$. A $\Q$-basis of $\Pic(S)$ is $\{\Lambda, N_\alpha\}$ where $\Lambda$ is the pullback of hyperplane class of $\kum(J(C))$ in $\PP^3$. Note that Keum wrote our $\Lambda$ as $H$. Then for each $\beta\in J(C)_2$:
\begin{align}
	\displaystyle T_\beta\sim \frac{1}{2}\left(\Lambda-\sum_{\alpha\in I(T_\beta)} N_\alpha\right).
\end{align}
The intersection products are given by $\Lambda^2=4, N_\alpha^2=-2$, and other products among $\Lambda, N_\alpha$ zero. One can deduce that $\Lambda\cdot T_\beta=2$, $T_\beta^2=-2$, and $T_\beta\cdot N_\alpha=1$ if $\alpha\in I(T_\beta)$, and $0$ otherwise.

To convert the above to the notation of \cite{Keum1997}, simply identify $p_6$ with $p_0$, and rename $T_{i6}$ by $T_i$ and $N_{i6}$ by $N_i$.
\subsection{Cones of divisors}
As a smooth K3 surface of Picard rank $17$, the birational geometry of $S$ is determined by its $(-2)$-curves. By a $(-2)$-curve we mean an irreducible curve $C\cong \PP^1$ on $S$. Such $C$ must have $C^2=-2$. There are infinitely many $(-2)$ curves on $S$. By \cite{Kovacs1994}, the Mori cone $\NEb(S)$ is spanned by all the $(-2)$-curve classes. The nef cone $\Nef(S)$ is dual to $\NEb(S)$. Furthermore, $\Nef(S)$ is locally polyhedral, with each $(-2)$-curve $C$ gives a codimension $1$ wall of $\Nef(S)$ with equation $D\cdot C=0$. No $(-2)$-curve is neglected in this correspondence:  by \cite[8.2.8]{Huybrechts2016}, if $H$ is ample, then $x:=H+(1/2)(H\cdot C)C$ is on the wall $C^\perp$, but not on any other walls of $C'\neq C$.
Now the Nakai-Moishezon criterion on $S$ reads:
\begin{align}
\amp(S)=\{D\in NS(S)\mid  D^2>0 \txt{and } D\cdot C>0 \txt{for all $(-2)$-curve } C\}.
\end{align}

\subsection{Classical involutions}\label{Cinv}
Let us review the generators of $\Aut(S)$. The details can be found in literatures \cite{Keum1997,Kondo1998,Dolgachev2002,Ohashi2009}.
\subsubsection{16 translations $t_\alpha$}
Each $\alpha\in J(C)_2$ induces a translation $x\mt x+\alpha$ on $J(C)$, which in turn gives the automorphism $t_\alpha$ of $S$. Each $t_\alpha$ is an involution, fixing $\Lambda$ and sending $N_\beta$ ($T_\beta$) to $N_{\beta+\alpha}$ ($T_{\beta+\alpha}$). The $16$ translations $t_\alpha$ all commute and form a group $(\Z/2\Z)^4$.

\subsubsection{The switch $\sigma$}
The Kummer surface $\kum(J(C))$ is self-dual and the projective isomorphism between $\kum(J(C))$ and its dual induces an involution $\sigma$ on $S$. This $\sigma$ exchanges each pair $(N_\alpha,T_\alpha)$, and sends $\Lambda$ to $3\Lambda-\sum_{J(C)_2}N_\alpha$.

\subsubsection{16 projections $p_\alpha$}
The projection of $\kum(J(C))$ from the node $N_\alpha$ gives a ramified double cover to $\PP^2$, which induces a covering involution $p_\alpha$ on $S$. The action of $p_\alpha$ on $NS(S)$ is just the reflection about the $(-4)$-root $\Lambda-2N_\alpha$.

\subsubsection{16 correlations $q_\alpha$}
One can construct $q_\alpha:=\sigma\circ p_\alpha \circ \sigma$, called the $16$ correlations. They are redundant as generators of $\Aut(S)$ but play a role when bounding a fundamental domain of $\Aut(S)$. The action of $q_\alpha$ on $NS(S)$ is the reflection about the $(-4)$-root $\sigma(\Lambda-2N_\alpha)$. 

\subsection{Hutchinson-G\"opel and Hutchinson-Weber involutions; Keum's automorphisms}

\subsubsection{G\"opel tetrads}
A G\"opel tetrad is a subset $g$ of $J(C)_2$ such that $\vv{g}=4$, and none of any three nodes in $g$ is contained in the same trope. Those tetrads in $J(C)_2$ which are not G\"opel are called the Rosenhain tetrads. The symmetric group $V=(\Z/2\Z)^4 \rtimes \mathcal{S}_6$ of the $(16,6)$-configuration is generated by the $16$ translations (the $(\Z/2\Z)^4$ part) and the symmetry group  $\mathcal{S}_6$ permuting the six nodes $N_0$ to $N_5$ on the same theta divisor. All the $60$ G\"opel tetrads form the same orbit under $V$.

Any G\"opel tetrad $g$ meets any $I(T_\alpha)$ at either $2$ or $0$ nodes. We say $g$ is of type 1 with respect to $T_\alpha$ if $\vv{g\cap I(T_\alpha)}=2$ and type 2 otherwise. In the following we always fix $T_0$ and drop the reference. Then all the G\"opel tetrads are:
\begin{align}
	\label{Gopel}
	g=\begin{cases}
		\{0,i6,jk,lm\} \txt{or } \{i6,j6,ik,jk\}, & \txt{type 1,}\\
		\{ik,il,jk,jl\}, & \txt{type 2.}\\
	\end{cases}
\end{align}
Clearly there are $45$ type 1 and $15$ type 2 G\"opel tetrads.

For each  G\"opel tetrad $g$ there is a Hutchinson-G\"opel involution $\varphi_g$ on $S$ \cite{Hutchinson1901}. \cite{Keum1997} referred to them by Cremona transformations, as they are induced by the Cremona of the  ambient $\PP^3$ which $\kum(J(C))$ embeds in via a G\"opel tetrad.

The action of $\varphi_g$ on $NS(S)$ is described in \cite[5.2]{Keum1997} \cite[5.1]{Kondo1998}. We have
\begin{lemma}
	Let $g$ be a type 1 G\"opel tetrad. Then $\varphi_g$ exchanges:
	\begin{align}
		\begin{cases}
			(T_0,T_{i6}), & \txt{if } g=\{0,i6,jk,lm\};\\
			(T_0,T_{ij}), & \txt{if } g=\{i6,j6,ik,jk\}.
		\end{cases}
	\end{align}
\end{lemma}
We list for example the action of $\varphi_g$ of type 1 on $NS(S)$ where $g=\{46,56,14,15\}$ for reference. The HG-involution $\varphi_g$ exchanges the following pairs:
\begin{align}
	(\Lambda, &\quad3\Lambda-2(N_{46}+N_{56}+N_{14}+N_{15}));\\
	(N_{\beta}, &\quad\Lambda-(N_{46}+N_{56}+N_{14}+N_{15})+N_{\beta}), \txt{for } \beta\in g.
\end{align}
\begin{alignat}{6}
	&&(N_{0}, \quad N_{23}),&&\quad (N_{16}, \quad N_{45});&& \quad (T_0, \quad T_{45}), &&\quad (T_{16}&, \quad T_{23});\\
	&&(N_{26}, \quad N_{12}),&&\quad (N_{36}, \quad N_{13});&&\quad (T_{46}, \quad T_{14}), &&\quad (T_{56}&, \quad T_{15});\\
	&&(N_{24}, \quad N_{25}),&&\quad (N_{34}, \quad N_{35});&&\quad (T_{24}, \quad T_{34}), &&\quad (T_{25}&, \quad T_{35}).
\end{alignat}
By linearity, we know the images of the remaining $4$ tropes:
\begin{align}
	\varphi_g(T_\beta) =\Lambda - (N_{46}+N_{56}+N_{14}+N_{15}) + T_\beta,  \txt{for } \beta\in g' = \{26,36,12,13\}.
\end{align}
The equations (9) and (13) hold in general, where 
\begin{align}
g':=
	\begin{cases}
		\{jl,jm,kl,km\}, & \txt{if } g=\{0,i6,jk,lm\};\\
		\{l6,m6,kl,km\}, & \txt{if } g=\{i6,j6,ik,jk\}.
	\end{cases}
\end{align}
\begin{defi}\label{zg}
	Let $g$ be a G\"opel tetrad of type 1. We define the Hutchinson-G\"{o}pel-type automorphism $z_g$ to be the automorphism $\varphi_g\circ t_{i6}$ if $g=\{0,i6,jk,lm\}$, or  $\varphi_g\circ t_{ij}$ if $g=\{i6,j6,ik,jk\}$.
\end{defi}

We note that $\varphi_g$ commutes with $t_{i6}$ or $t_{ij}$ in the respective case, so $z_g$ itself is an involution. By definition $z_g$ fixes $T_0$. With the action of $\varphi_g$ above we find:
\begin{lemma}
	\label{zgexample}
	Let $g=\{46,56,14,15\}$. Then $z_g$ fixes the following classes:
	\[ T_0, T_{16}, T_{23}, T_{45}, N_{24},N_{34}, N_{25},N_{35},\]
	and exchanges the following pairs:
\begin{align}
	(\Lambda, &\quad3\Lambda-2(N_{46}+N_{56}+N_{14}+N_{15}));\\
	(N_{\beta}, &\quad\Lambda-(N_{46}+N_{56}+N_{14}+N_{15})+N_{\beta+(45)}), \txt{for } \beta\in g;\\
	(T_{\beta}, &\quad\Lambda-(N_{46}+N_{56}+N_{14}+N_{15})+T_{\beta+(45)}), \txt{for } \beta\in g' = \{26,36,12,13\}.
\end{align}
\begin{alignat}{6}
	&&(N_{0}, \quad N_{16}),&&\quad (N_{26}, \quad N_{36});&& \quad (T_{46}, \quad T_{15}), &&\quad (T_{56}&, \quad T_{14});\\
	&&(N_{12}, \quad N_{13}),&&\quad (N_{23}, \quad N_{45});&&\quad (T_{24}, \quad T_{35}), &&\quad (T_{25}&, \quad T_{34}).
\end{alignat}
\end{lemma}

\subsubsection{Weber hexads}
A Weber hexad is a subset $w$ of $J(C)_2$ such that $\vv{w}=6$ and no four of the six nodes are in a trope or a G\"opel tetrad. Equivalently, $w$ is the symmetric difference of a G\"opel and a Rosenhain tetrad. There are in total $192$ Weber Hexads, which form $12$ orbits under translations.

A Weber hexad $w$ meets $I(T_\alpha)$ at either $3$ or $1$ nodes. We say $w$ is of type 1 with respect to $T_\alpha$ if $\vv{w\cap I(T_\alpha)}=3$ and type 2 otherwise. Same as the G\"opel tetrads, we will fix $T_0$ as the reference. Then all the Weber tetrads are classified by:
\begin{align}
	\label{Weber}
	w=\begin{cases}
		\{0,i6,j6,ik,kl,lj\}, \txt{or } \{i6,j6,k6,ij,il,jm\}, & \txt{type 1,}\\
		\{0,ij,jk,kl,lm,mi\}, \txt{or } \{i6,ij,jk,ki,jl,km\} & \txt{type 2.}\\
	\end{cases}
\end{align}
There are then $120$ Weber hexads of type 1 and $72$ of type 2. Each $w$ of type 1 has a dual $w'$, which is the unique Weber hexad, which is also of type 1, such that the symmetric difference $w\triangle w'=I(T_0)$. Specifically, $\{0,i6,j6,ik,kl,lj\}$ and $\{k6,l6,m6,ik,kl,lj\}$ are dual to each other. The $120$ Weber hexads form $60$ pairs of dual hexads in this sense.

To each Weber hexad $w$ there is a Hutchinson-Weber involution (HW for short) $\varphi_w$. Their actions on $NS(S)$ are very complicated. The description that $\varphi_w$ exchanges ten nodes with ten tropes can be found in \cite{Dolgachev2002} as Enrique's involutions on a Hessian surface. See \cite[\S 7, Remark]{Ohashi2009} for a complete description. 

Keum \cite{Keum1997} first constructed some infinite order automorphisms $k_w$ for each $w$. We introduce a different definition when $w$ is of type $1$.
\begin{defi}
	For each Weber hexad $w$ of type 1, we define the Keum's automorphism $z_w$ to be the automorphism $\varphi_w\circ t_\alpha \circ \sigma=\varphi_w\circ \sigma \circ t_\alpha$, where $\alpha=ij$ if $w=\{0,i6,j6,ik,kl,lj\}$, or $lm$ if $w=\{i6,j6,k6,ij,il,jm\}$.
\end{defi}

Some tedious calculation shows that $\varphi_w(T_0)=N_{ij}$ when $w=\{0,i6,j6,ik,kl,lj\}$, and $N_{lm}$ when $w=\{i6,j6,k6,ij,il,jm\}$. Therefore $z_w(T_0)=T_0$ in either case. Furthermore, the action of $z_w$ coincides with Keum's $k_w$. Hence by the global Torelli Theorem on K3 surfaces, they are exactly the automorphisms in \cite{HY2019} which lift to pseudo-automorphisms of the threefold $X$.

\begin{theorem}\cite{Kondo1998} The automorphism group of a Jacobian Kummer surface $S$ of Picard number $17$ is generated by the $16$ translations $t_\alpha$, the switch $\sigma$, the $16$ projections $p_\alpha$, the $60$ Hutchinson-G\"opel involutions $\varphi_g$, and the $192$ Keum's automorphisms $k_w$.
\end{theorem}
\begin{remark}
\label{hgw relation} 
In \cite{Kondo1998} the author introduced a different definition of Keum's $k_w$, based on the theory of fundamental domains constructed out of Leech roots.
Each Leech root $r$ is associated with two Weber Hexads $\Phi_r$ and $\Psi_r$. The hexads induce an automorphism $\psi_r$. If $w$ is of type 1, then those $z_w$ that we defined coincides with those $\psi_r$: let $w = \Psi_r$ be of type 1, then $z_w = \psi_r = k_w$. This gives a new connection between the Hutchinson-Weber involutions $\phi_w$ and Kond\=o's construction of $\psi_r$. To prove this, compare their actions on the $\mathbb{Q}$-basis $\{H,N_\alpha\}$. Since $\phi_w$ exchanges $10$ nodes with $10$ tropes, composing with the switch $\sigma$ gives $z_w$ which maps the $10$ nodes to other $10$ nodes. This is then easily compared with the description of the isometry by $\psi_r$ in  \cite[6.2]{Kondo1998}.
\end{remark}
\section{HG-type Pseudo-automorphisms}
\subsection{The threefold and the anticanonical surface}
In this section, we construct the HG-type pseudo-automorphism $\phi$ out of the underlying involution.
Let $X$ be the successive blow-up of $\PP^3$ along $6$ very general points $p_1,\cdots, p_6=p_0$ and $15$ lines $l_{ij}:=\overline{p_i p_j}$. The Jacobian Kummer $S$ embeds in $X$ as the unique anticanonical section. Let $R$ be the rational normal curve in $X$ through all $p_i$. Then $R=T_0$. The restriction to $S$ induced the map $r:\Pic(X)\ra\Pic(S)$, given by \cite[Prop 5.1]{HY2019}:
\begin{alignat}{4}
	H & \mt \rlap{$(3\Lambda- \sum_{1\leq i<j\leq 5} N_{ij})/2,$} && &&\\
	E_i & \mt N_{i6}, \quad && E_0\mt N_0,\quad  && E_{ij}  \mt T_{ij}.
\end{alignat}

\begin{prop}\label{injhom}
	The group homomorphism $u:PsAut(X)\ra \Aut(S)$ is injective.
\end{prop}
\pf. We prove it in 3 steps. Suppose $u(f)=\id$. 

\underline{Step 1.} We show $f(E_i)=E_i$ and $f(E_{ij})=E_{ij}$. Indeed, suppose $C=f(E_i)$. Then $C$ is an irreducible effective divisor in $X$. Restricting to $S$, we find $r(C)=E_i$. Therefore $r(C-E_i)\sim 0$. Now by \cite[Lemma 5.6]{HY2019}, $r(C-E_i)\sim 0$ implies that $\deg (C-E_i)=0$, so $\deg C=0$. Now the only degree $0$ irreducible effective divisors on $X$ are $E_i$ and $E_{ij}$, and $r(E_{ij})=T_{ij}\neq E_i$. Hence $C=E_i$.
 Similar argument shows that $f(E_{ij})=E_{ij}$.

\underline{Step 2.} Now we conclude that $f(H)=H$. Indeed, this follows from Step 1 and that $f(K_X)=K_X$.

\underline{Step 3.} By Step 2, the action of $f$ on $\Pic(X)$ is identity. Consider the birational map $g:\PP^3\dashrightarrow \PP^3$ induced by $f$. Since $f(H)=H$, we have $g$ is linear, so $g\in \PGL(4)$. Since $5$ general points decides an element in $\PGL(4)$, and $g$ fixes the $6$ points $p_0$ to $p_5$, we find that $g$ is the identity map. Hence $f=\id$ and $u$ is injective.\qed

\begin{theorem}
	\label{mainG}
	For each G\"opel tetrad $g$ of type 1, there is a pseudo-automorphism $\phi_g\in PsAut(X)$ such that $u(\phi_g)=z_g$ and $\phi_g$ itself is an involution. Each $\phi_g$ is induced by the complete linear system 
	\begin{align}
D_g:=5H-2\sum_{i=0}^5 E_i-2\sum_{\alpha\in g'} E_\alpha,
	\end{align}
where 
\begin{align}
g':=
	\begin{cases}
		\{jl,jm,kl,km\}, & \txt{if } g=\{0,i6,jk,lm\};\\
		\{0l,0m,kl,km\}, & \txt{if } g=\{i6,j6,ik,jk\}.
	\end{cases}
\end{align}
\end{theorem}
Clearly all the $45$ $\phi_g$ are symmetric under the permutation of the six points $p_i$. In fact, a careful calculation shows that the actions of those $z_g$ (not those $\varphi_g$) on $NS(S)$ are also symmetric under the automorphisms of the (16,6)-configuration that relabel the six nodes on $T_0$ (see \cite[1.7]{Keum1997}). As a result, in the following we fix $g=\{46,56,14,15\}$, with $g'=\{02,03,12,13\}$, so $D=D_g=5H-2\sum_i E_i-2(E_{02}+E_{03}+E_{12}+E_{13})$, and reduce the proof of Theorem \ref{mainG} to this very case.

\subsection{The complete linear system}
Write $\phi:=\phi_g$ for short. We first claim the action of $\phi$ on $\Pic(X)$:
\begin{prop}
	\label{PicXg}
	The involution $\phi$ acts on $\Pic(X)$ by exchanging the following pairs:
\begin{alignat*}{7}
	&& (H, &\quad D), && \quad (E_{0}, &\quad E_1), && \quad (E_2, &\quad E_3),\\
	&& (E_{04}, &\quad E_{15}), && \quad (E_{05}, &\quad E_{14}), &&\quad  (E_{24}, &\quad E_{35}), && \quad (E_{34}, &\quad E_{25});
\end{alignat*}
\begin{align*}
	(E_{4}, &\quad F_5:=2H-\sum\nolimits_{i=0}^5 E_i+E_5 -\sum\nolimits_{\alpha\in g'}E_\alpha), \\
	(E_{5}, &\quad F_4:=2H-\sum\nolimits_{i=0}^5 E_i+E_4 -\sum\nolimits_{\alpha\in g'}E_\alpha ), \\
	(E_{\alpha}, &\quad H_\beta:=2H-\sum\nolimits_{i=0}^5 E_i-\sum\nolimits_{\alpha\in g'}E_\alpha +E_{\beta}), \txt{where }\alpha, \beta \in g', \alpha\cap \beta=\emptyset.
\end{align*}
and fixing the following divisors:
\[E_{01},E_{23}, E_{45}.\]
\end{prop}

We prove Theorem \ref{mainG} following the same strategy in \cite{HY2019}. The steps are as follows.
\begin{enumerate}
	\item Show $h^0(D)=4$ and find various generating sets of sections of $\vv{D}$.
	\item Show there are six points $q_i$ in $\PP^3$ such that the map $\phi_D:X\dashrightarrow \PP^3$ maps $R$ to the unique rational normal curve $R'$ through $q_i$. Deduce that $\{p_i\}$ and $\{q_i\}$ are projectively equivalent.
	\item Show $\phi_D\circ \phi_D=\id $ on an open subset in $\PP^3$. Then $\phi_D$ is birational and is an involution.
	\item Compute the Jacobian determinant of $\phi_D$ to show it contracts exactly the quadrics $F_4,F_5,H_{\beta}$ above to points and lines.
	\item After blowing-up the points and lines on the target, show $F_5$ ($F_4$) is mapped to $E_4$ ($E_5$), and $H_{\beta}$ to $E_{\alpha}$. Show $F_5$ and $H_{\beta}$ are not contracted.
	\item Finally show the remaining $11$ lines are either birationally paired or fixed.
\end{enumerate}
In particular we will apply the key results  \cite[8.5, 8.9, 9.2]{HY2019}.
\subsection{The sections and birationality}
\begin{lemma}
	\label{fhunique}
	For general six points $p_i$, the quadric classes $F_i$ and $H_{ij}$ have unique irreducible and distinct sections. Subtracting any $E_k$ or $E_{kl}$ from $F_i$ or $H_{ij}$ results in a divisor not effective.
\end{lemma}
\pf. Direct calculation on their sections.\qed

\begin{lemma}
	\label{resFH}
	The action table in Proposition \ref{PicXg} restricts to the action of $z_g$ on $S$. In particular,
	\begin{align*}
	&D_{\mid S}\sim z_g(H_S), \\
	&{F_4}_{\mid S}\sim \Lambda-(N_{46}+N_{14}+N_{15}),\\
	&{F_5}_{\mid S}\sim \Lambda-(N_{56}+N_{14}+N_{15}),\\
	&{H_\beta}_{\mid S}\sim z_g(T_\alpha).
	\end{align*}
\end{lemma}
\pf. Direct calculation.\qed

We introduce notations. Up to a scalar, let $p_{ijk}$ be the polynomial defining the plane $\Gamma_{ijk}$. Write $f_i$ and $h_{ij}$ for the unique section of $F_i$ and $ H_{ij}$. Write $x_i$ or $x_{ij}$ for the unique section of the exceptional divisor $E_i$ or $E_{ij}$.
\begin{lemma}
\label{secD}
	For very general six points $p_i$, $h^0(X,D)=4$. $\vv{D}$ can be generated by the following sections:
	\begin{equation}
		\label{s01}
			(p_{023} h_{02} h_{03}x_0 x_2 x_3 x_{23}, p_{123} h_{12} h_{13}x_1 x_2 x_3 x_{23},  p_{012} h_{02} h_{12}x_0 x_1 x_2 x_{01}, p_{013} h_{03} h_{13}x_0 x_1 x_3 x_{01}),
	\end{equation}
	or
	\begin{equation}
		\label{s02}
(p_{035} f_5 h_{03}x_0 x_3 x_{05}x_{35}, p_{025} f_5 h_{02}x_0 x_2 x_{05}x_{25},  p_{134} f_4 h_{13}x_1 x_3 x_{14}x_{34}, p_{124} f_4 h_{12}x_1 x_2 x_{14}x_{24}),
	\end{equation}
	or
	\begin{equation}
		\label{s03}
(p_{035} f_5 h_{03}x_0 x_3 x_{05}x_{35}, p_{025} f_5 h_{02}x_0 x_2 x_{05}x_{25},  p_{125} f_5 h_{12}x_1 x_2 x_{15}x_{25}, p_{124} f_4 h_{12}x_1 x_2 x_{14}x_{24}).
	\end{equation}

\end{lemma}

\pf. Direct calculation shows that each term is a section of $D$. The four sections  in (\ref{s01}) are defined on an open set in $\PP^3$ by the polynomials
\[(s_0,s_1,s_2,s_3):=(p_{023} h_{02} h_{03}, p_{123} h_{12} h_{13},  p_{012} h_{02} h_{12}, p_{013} h_{03} h_{13}).\]
To show the linear independence, we need only prove that $s_i$ are linearly independent. Suppose there are scalars $a_i$ such that $\sum_i a_i s_i=0$. Say $a_0\neq 0$. Then $h_{02}\mid a_1  p_{123} h_{12} h_{13}+   a_3 p_{013} h_{03} h_{13}$. Since $h_{02}$ are $h_{13}$ are not scalar multiple to each other,
\[h_{02} \mid a_1  p_{123}  h_{12}+  a_3 p_{013} h_{03}.\]
Therefore there exists a linear term $\beta$ such that $\beta h_{02} h_{13} = a_1  p_{123}  h_{12}+  a_3 p_{013} h_{03}$. Examine now the orders of vanishing at $p_i$ and $l_{ij}$ of both side, and note that we can choose the open set in the previous step to include those $p_i$ and $l_{ij}$ different from the exceptional classes appearing in $s_i$. It shows that $\beta$ must vanish at both the lines $l_{02}$ and $l_{13}$, which is impossible. Hence $a_0=0$. By symmetry all $a_i=0$.

Now we know $h^0(X,D)\geq 4$. On the other hand, we use the short exact sequence:
\[0\ra H^0(X, D-S)\ra H^0(X,D)\ra H^0(S,D_{\mid S}).\]
Clearly $D-S$ is not effective, so $h^0(X,D)\leq h^0(S,D_{\mid S})$. Now Lemma \ref{resFH} says $D_{\mid S}\sim z_g(H_S)$.
Since $h^0(S,H_S)=4$ \cite[Prop. 5.7]{HY2019} and  $z_g$ is an automorphism, we know $h^0(S,D_{\mid S})=4$. Hence $h^0(X,D)=4$.

Finally (\ref{s02}) and (\ref{s03}) are linearly independent by a similar argument.\qed

\begin{lemma}
	\label{Dmov}
	The class $D$ is movable.
\end{lemma}
\pf. From (\ref{s01}), the base locus of $\vv{D}$ lies in the pairwise intersections among the divisors $E_i, E_{jk}, \Gamma_{ijk}$ and $H_\beta$, which do not contain any divisor. \qed

Now consider $\phi_D:X \dashrightarrow \PP^3$. We show
\begin{lemma}
	\label{4points}
	The map $\phi_D$ contracts $E_i$, $i=0,1,2,3$ to $4$ distinct points $q_j$, $j=1,0,3,2$ respectively, and contracts each $H_{ij}$ into the line $\overline{q_k q_l}$ where $\{i,j,k,l\}=\{0,1,2,3\}$.
	Let $\psi:\PP^3\dashrightarrow\PP^3$ be the map induced by $\phi_D$, where $p_i$ is identified with $q_i$. Then $\psi$ maps the line $l_{01}$ birationally to the line $\overline{q_0 q_1}$.
\end{lemma}
\pf. We read directly in (\ref{s01}) that
\begin{align*}
	E_0\mt &\quad [0:1:0:0]=:q_1,&	E_1\mt &\quad [1:0:0:0]=:q_0,\\
	E_2\mt &\quad [0:0:0:1]=:q_3,&	E_3\mt &\quad [0:0:1:0]=:q_2.
\end{align*}
Then $H_{02}$ is contracted to some locus in the line $\overline{q_1 q_3}$. The other $H_{ij}$ (where $ij \in g'$) follow from symmetry.
Finally, a local calculation restricting $[s_0:s_1:s_2:s_3]$ on the line $l_{01}$ shows that $l_{01}$ is mapped isomorphically to its image $\overline{q_0 q_1}$. 
\qed

\begin{lemma}
	\label{2points}
	The map $\phi_D$ contracts $F_4$ and $F_5$ to two points $q_5$ and $q_4$ respectively. If the six points $p_i$ are very general, then $\{p_i\}$ and $\{q_i\}$ are projectively equivalent. In particular $\{q_i\}$ are all distinct.
\end{lemma}
\pf. Use the sections (\ref{s03}). Then we see $F_5$ is contracted to a point $q_4$, and by symmetry $F_4$ is contracted to a point $q_5$. 

We show $\phi_D(R)=R'$ is a rational normal curve through the six points $q_i$. Similar to \cite[Thm 8.4]{HY2019}, we set $p_i$ at 
\begin{alignat}{5}
	    p_0=[1:0:0:0],  &&\quad  p_1=[0:1:0:0],  &&\quad p_2=[0:0:1:0],\\
	    p_3=[0:0:0:1], &&\quad p_4=[1:1:1:1], &&\quad p_5=\left[1:\frac{1}{a}:\frac{1}{b}:\frac{1}{c}\right].
\end{alignat}
for $a,b,c$ nonzero. Then the rational normal curve $R_0$ in $\PP^3$ through $p_i$ is parameterized by 
\begin{equation}
	[u:v]\mt \left[\frac{1}{u+v}: \frac{1}{au+v}:\frac{1}{bu+v}:\frac{1}{cu+v}\right].
	\label{Rpara}
\end{equation}
Identifying $R_0$ with $R$, we find the $[u:v]$ coordinates of $p_i$ are respectively
\begin{align}
	u/v=(-1,-1/a,-1/b,-1/c,0,\infty).
\end{align}
Now the HG-type involution $z_g$ sends $(N_0,N_{16},\cdots, N_{56})$ to \[U:=(N_{16},N_0,N_{36},N_{26},z_g(N_{46}),z_g(N_{56})).\] Since $z_g(R)=R$, we find the six points $\{p_{i}\}$, as the unique intersection of $R$ with $\{N_{i6}\}$, are projectively equivalent to their images under $z_g$, namely: $\{p_1,p_0,p_3,p_2,x,y\}$, where 
\[\{x\}:=z_g(N_{46})\cap R=F_5\cap R, \quad \{y\}:=z_g(N_{56})\cap R=F_4\cap R,\]
by Lemma \ref{resFH}.

Now we can compute the image $R'$, by plugging (\ref{Rpara}) into the polynomial of the sections of $\phi_D$. Denote by $p_{ijk,R}$ and $h_{ij,R}$ the polynomials obtained by this process. Then these polynomials are computed by counting the multiplicities at the six points $p_i$. We have:
\begin{align}
& h_{02,R} = h_{03,R} = h_{13,R} = h_{23,R} = uv(u+v)(au+v)(bu+v)(cu+v),\\
& p_{023,R} = (u+v)(bu+v)(cu+v), \quad  p_{123,R} = (au+v)(bu+v)(cu+v),  \\
& p_{012,R} = (u+v)(au+v)(bu+v), \quad  p_{013,R} = (u+v)(au+v)(cu+v).
\end{align}

As a result, the curve $R'$ is given by
\begin{align}
	[u:v]\mt \left[\frac{1}{au+v}: \frac{1}{u+v}:\frac{1}{cu+v}:\frac{1}{bu+v}\right].
\end{align}
Clearly $q_0,\cdots,q_3\in R'$. Now $\phi_D(F_5)=q_4$. On the other hand, we show that $\phi_D$ is defined at $x$, so $\phi_D(x)=q_4$. By symmetry $\phi_D(y)=q_5$. With this, we will find that the $[u:v]$ coordinates of $q_i$ equal to the $[u:v]$ coordinates of $\{p_1,p_0,p_3,p_2,x,y\}$ respectively. Hence $q_i$ are projectively equivalent to $p_i$ and all $q_i$ are distinct, which completes the proof.

Let us show that $\phi_D$ is defined at $x$. Indeed, we need only show that at least one of the sections in (\ref{s01}) does not vanish at $x$. Consider first those $h_{ij}$ restricted to $R$. By definition, $h_{ij}$ vanishes at $z \in R$ if and only if $\{z\}=H_{ij}\cap R = z_g(T_{kl} \cap R)$, where $\{i,j,k,l\} = \{0,1,2,3\}$. But any two tropes are disjoint on $S$, and $R=T_0$, so all $h_{ij}$ are non-vanishing on $R$. Next, at least one of the planes in (\ref{s01}) is away from $x$. Finally, among the exceptional sections $x_i$ and $x_{ij}$, $x_{ij}$ does not vanish on $R$ by the same argument above, and that $E_{ij}\mid_{S} = T_{ij}$. For $x_i$ where $i\in\{0,1,2,3\}$, $E_i\mid_{S} = N_{i6}$ meets $R$ at a different point from $x$, which is where $F_5$ meets $R$, because under the involution $z_w$, $N_{46}$ does not meet $N_{i6}$ for $i\in\{0,1,2,3\}$. This completes the argument that at least one section in (\ref{s01}) does not vanish at $x$. 
\qed

\begin{prop}\label{selfComposition}
	The composition $\phi_D\circ \phi_D$ is identity on an open set in $\PP^3$. Thus $\phi_D$ is a birational involution. The exceptional set of $\phi_D$ equals the union of the six quadrics $F_4,F_5$, and $H_{ij}$ for $ij \in g'$.
\end{prop}
\pf. We need only compute the polynomials $h_{ij}(\phi_D)$ and $p_{ijk}(\phi_D)$. By \cite[Prop. 8.5]{HY2019}, we can read the multiplicities of $h_{ij}$ and $f_i$ in them by counting the multiplicities of the corresponding lines or points in the divisor classes of $H_{ij}$ or $\Gamma_{ijk}$. For instance, we have $f_4,f_5,h_{02},h_{03}$ and $h_{12}$ all divide $h_{02}(\phi_D)$. Here $h_{02}(\phi_D)$ is a nonzero polynomial since by Lemma \ref{4points}, $\phi_D(E_{01})=\overline{q_0 q_1}$, which is not contained in $H_{02}$ on the target. Then by degree comparison, up to a nonzero scalar:
\[h_{02}(\phi_D)=f_4 f_5 h_{02} h_{03}h_{12}.\]
Let $A:=f_4 f_5  h_{12} h_{13}h_{02}h_{03}$. Then $A=h_{02}(\phi_D)h_{13}$. 
By symmetry, $A=h_{03}(\phi_D)h_{12}=h_{12}(\phi_D)h_{03}=h_{13}(\phi_D)h_{02}$. Similarly we find
\begin{alignat*}{2}
	& p_{023}(\phi_D)=h_{12} h_{13} p_{123},\quad && p_{123}(\phi_D)=h_{02} h_{03} p_{023},\\
	& p_{012}(\phi_D)=h_{03} h_{13} p_{013},\quad &&p_{013}(\phi_D)=h_{02} h_{12} p_{012}.
\end{alignat*}
Therefore when we multiply these polynomials, we find on an open locus of $\PP^3$,
\[\phi_D\circ \phi_D=[p_{123}:p_{023}:p_{013}:p_{012}]=[x:y:z:w]\]
if we place $(p_0,\cdots,p_3)$ at $([1:0:0:0],\cdots, [0:0:0:1])$. Therefore $\phi_D\circ\phi_D$ is the identity map. This shows that $\phi_D$ is birational.

	Finally we can determine the exceptional set of $\phi_D$ via computing the Jacobian determinant $J$ of $\phi_D$. From Lemma 8.9 of \cite{HY2019}, $f_4^2,f_5^2$ and each $h_{ij}$ divides $J$. Since $\phi_D$ is birational, $J\not\equiv 0$. Now $J$ has degree $(5-1)\times 4=16$, so by a degree comparison: 
	\[J=f_4^2 f_5^2 h_{02}h_{03}h_{12}h_{13}.\]
	Hence the exceptional set of $\phi_D$ contains only these quadrics.\qed
\subsection{The exceptional loci}
We have proved the steps (1) to (4) and here we show (5) and (6).
\begin{lemma}
	\label{f4h02}
	\begin{enumerate}
		\item	After blowing-up the point $q_4$, the induced birational map $\phi_4:X\dashrightarrow \bl_{q_4} \PP^3$ maps $F_5$  to $E_4$, and does not contract $F_5$.
		\item	After blowing-up the line $\overline{q_1 q_3}$, the induced birational map $\phi_{13}:X\dashrightarrow \bl_{\overline{q_1 q_3}} \PP^3$ maps $H_{02}$ to $E_{13}$, and does not contract $H_{02}$.

	\end{enumerate}
\end{lemma}
\pf. We apply \cite[9.1, 9.2]{HY2019}. All we need is a local calculation at $F_5$ and $H_{02}$, and compute the Jacobian determinant on an affine chart. Then the method is the same as the proof of \cite[9.3]{HY2019}.

(1). Using the sections (\ref{s02}), and the local map $\xi:U\dashrightarrow E_4$, we find  up to a nonzero scalar:
\begin{align}
	\det J(\xi)=\frac{x_0}{(s_0 s_1)^2}\det J(s_0,s_1,s_2,s_3)_{x_0,x_1,x_2,x_3}= \frac{x_0  h_{02}h_{13} f_4^2 }{ p_{035} p_{025}},
\end{align}
which does not vanish at the generic point of $F_5$. Hence $F_5$ is not contracted.

(2). Using the sections (\ref{s01}), and the local map $\xi:V\dashrightarrow E_{13}$, we find
\begin{align}
	\det J(\xi)=\frac{x_0}{s_0 s_1^3}\det J(s_0,s_1,s_2,s_3)_{x_0,x_1,x_2,x_3}= \frac{x_0  h_{12} f_5^2 f_4^2 }{ p_{023} p_{013}^3 h_{13}^2 h_{03}^3},
\end{align}
which does not vanish at the generic point of $H_{02}$. Hence $H_{02}$ is not contracted.\qed

\begin{lemma}
	\label{e13}
After blowing-up the point $q_1$, the induced birational map $\phi_1:X\dashrightarrow \bl_{q_1} \PP^3$ maps $E_0$  to $E_1$, and does not contract $E_0$.
\end{lemma}
\pf. This is a local calculation at $p_0$, and we need only blow up the lines appearing in $D$, namely $l_{02}$ and $l_{03}$. Then $E_0$ is isomorphic to the blow-up of $\PP^2$ at two points $a_2$ and $a_3$, so $\Pic(E_4)=\Z\{h, e_2,e_3\}$, with $e_i$ the exceptional divisor over $a_i$. Let $\ell_{ij}$ be the intersection of the proper transform of $l_{ij}$ with $E_0$. Then Table \ref{e4r} shows the restrictions of the defining polynomials we will use.
\begin{table}[ht]
    \centering
\[\begin{array}{r|l|l}
		\xi\in \Pic(X)	 & r_4(\xi) & \txt{the zeroes of the restricted section} \\
		\tilde{\Gamma}_{023} & h-e_2-e_3  & \ell_{23}\\
		\tilde{\Gamma}_{012} & h-e_2  & \ell_{12}\\
		\tilde{\Gamma}_{013} & h-e_3  & \ell_{13}\\
		H_{02} & h-e_3 & c\\
		H_{03} & h-e_2 & d\\
		H_{12} & h-e_2-e_3 & \ell_{23}\\
		H_{13} & h-e_2-e_3 & \ell_{23}\\
\end{array}\]\caption{Restrictions to $E_0$}
    \label{e4r}
\end{table}
Here $c$ and $d$ are linearly independent linear terms. Then up to nonzero scalars, the induced map $\sigma:E_0\dashrightarrow E_2$ is given by
\begin{align}
	[p_{023} h_{02} h_{03}: p_{012} h_{02} h_{12} : p_{013} h_{03} h_{13}]_{\mid E_0}&=[cd:c\ell_{12}:d\ell_{13}],
\end{align}
which is the standard Cremona map $\PP^2\dashrightarrow \PP^2$ after a change of coordinates $[X:c:d]\mt [X:Y:Z]$. Hence $\sigma$ is birational onto $E_1$.\qed

\begin{lemma}
	\label{11lines}
	After blowing-up all the $15$ lines on the target, the induced map $\phi:X\dashrightarrow X$ fixes $E_{01},E_{23}, E_{45}$, and permutes the pairs
	\[ (E_{04},  E_{15}),  (E_{05},  E_{14}),   (E_{24},  E_{35}),  (E_{34},  E_{25}).\]
\end{lemma}
\pf. By Lemma \ref{selfComposition}, none of the $11$ lines above are in the exceptional set of $\phi_D$, and hence are not contracted by $\phi_D$. By Lemma \ref{4points}, $\phi_D$ fixes $l_{01}$ and $l_{23}$. Since $\phi_D$ is an involution, it must be that $\phi(E_{01}) = E_{01}$ and  $\phi(E_{23}) = E_{23}$. 

Let $\psi:\PP^3\dashrightarrow\PP^3$ be the induced map. For the remaining lines, we need only show $\psi$ permute the pairs of the underlying lines, and fixes $l_{45}$, by the same reason above.

We first prove that $\psi(\Gamma_{034})=\Gamma_{125}$. Indeed, we need show that $p_{125}(\psi)=p_{034}f_4 h_{03}$. Then either by symmetry or that $\psi$ does not contract $\Gamma_{034}$, we have the equality.

Since $\psi$ is birational, $p_{125}(\psi)$ is not a zero polynomial. Since $\psi$ contracts $F_4$ to $q_5$ and $H_{03}$ to $\overline{q_1 q_2}$, we have 
\[p_{125}(\psi)=f_4 h_{03} \beta\]
for some linear term $\beta$.
Compare the vanishing multiplicities along the points. Then $\beta$ vanishes at $E_4$ and $E_{03}$. Hence $\beta=p_{034}$. This proves the claim.

Now we exploit symmetry so that:
\begin{align}
	\psi(\Gamma_{034})=\Gamma_{125}, &  \quad \psi(\Gamma_{134})=\Gamma_{025},\\
	\psi(\Gamma_{024})=\Gamma_{135}, &  \quad \psi(\Gamma_{124})=\Gamma_{035}.
\end{align}
Hence $\psi(l_{14})\subset \Gamma_{025}\cap \Gamma_{035}=l_{05}$. Symmetry shows $\psi(l_{14})=l_{05}$ and the other three pairs are exchanged similarly.
Finally, $p_{145}(\psi)$ is divided by $f_4 f_5$. So $p_{145}(\psi)=\eta f_4 f_5$ where $\eta$ is a linear term. Similar argument shows that $\eta$ vanishes at $p_4$ and $p_5$. Furthermore, $p_{145}(\psi(l_{04}))=p_{145}(l_{15})=0$. Therefore $\eta=p_{045}$, as neither $F_4$ nor $F_5$ passes throught the line $l_{04}$. Hence $\psi(\Gamma_{045})=\Gamma_{145}$, so by symmetry,  $\psi(\Gamma_{245})=\Gamma_{345}$, and hence $\psi$ fixes $l_{45}$.\qed

\pfof{Theorem \ref{mainG}}. 
As a summary, we proved that $\phi$ is a birational automorphism of $X$ and $\phi$ contracts no divisors. Then by \cite[2.1]{Cantat2017}, $\phi$ is a pseudo-automorphism. Proposition \ref{PicXg} then follows from the proofs above and Lemma \ref{fhunique}, \ref{Dmov}. Now by Lemma \ref{resFH}, $u(\phi)$ has the same action as $z_g$ on $NS(S)$. The global Torelli theorem of K3 surfaces implies that $u(\phi)=z_g$. This finishes the proof of Theorem \ref{mainG}.\qed

\section{The group of HG-type and Keum's automorphisms}
In this section we prove Theorem \ref{Gstab} following Kond\=o's proof \cite{Kondo1998} of generators of $\Aut(S)$.

\subsection{Leech lattice and Leech roots}
This paragraph reviews Kond\=o's construction. There is a unique unimodular even lattice $\II_{1,25}$ of signature $(1,25)$. A root of $\II_{1,25}$ is an element $r$ with $\vr{r,r}=-2$. The roots in $\II_{1,25}$ are one-to-one corresponding to the elements in the (negative definite) {\em Leech lattice} $L$ by Conway. 
Let $U$ be the hyperbolic lattice. Then $\II_{1,25}\cong L\oplus U$. The Weyl vector $w=(0,(0,1))$ with $0\in L, (0,1)\in U$, is a vector in $\II_{1,25}$.
A {\em Leech root} of $\II_{1,25}$ is a root $r$ such that $\vr{r,w}=1$.
By \cite{Conway1983}, the Weyl group $W(\II_{1,25})^{(2)}$ of reflections about the Leech roots acts on the positive cone of $\II_{1,25}$ with a fundamental domain $D$:
\begin{align}
	D=\{x\in \II_{1,25}\mid \vr{x,x}>0 \txt{and } \vr{x,r}>0 \txt{for all  Leech roots } r\}.
\end{align}
This $D$ is rational polyhedral, with infinitely many faces. The isometry group $O(\II_{1,25})$ 
is the split extension of $W(\II_{1,25})^{(2)}$ by $Sym(D)$, the group of isometries which keeps $D$ invariant.

The connection between the Leech roots to Kummer surfaces lie in their role in constructing a fundamental domain of $\Aut(S)$ acting on the ample cone of $S$. The Picard lattice $NS(S)$ has signature $(1,16)$, and can be primitively embedded into the lattice $\II_{1,25}$, such that $NS(S)^\perp$ is generated by Leech roots. By \cite{Kondo1998}, a fundamental domain is 
\[D':=D\cap P(S),\]
where $P(S)$ is the positive cone of $S$, i.e., the branch of the cone $\{x\in NS(S)\mid x\cdot x>0\}$ with an ample class.
In fact, $D'\subset \overline{\amp(S)}$, since $D'$ contains the projection $w'$ of the Weyl vector $w$, and $w'$ is an ample divisor:
\begin{align}
	w'=\frac{1}{4}\left(\sum_\alpha N_\alpha +\sum_\alpha T_\alpha\right)\sim 2\Lambda-\frac{1}{2}\sum_\alpha N_\alpha,
\end{align}
and then every $f\in\Aut(S)$ lifts to an isometry in $O(\II_{1,25})$, which either fixes $D$ or send $D$ to another fundamental domain.

Next, Borcherds' \cite{Borcherds1987} methods shows that $D'$ is rational polyhedral, with only finitely many faces, each bounded by a Leech root $r$ such that $r^\perp$ in $\II_{1,25}$ does not miss $NS(S)$.

In general for a projective $K3$ surface $S$, we have a group homomorphism:
\begin{align}
	\Aut(S)\ra O(NS(S))/W(NS(S))^{(2)}\cong Sym (\amp(S)).
\end{align}
with kernel finite.
Here for the Kummer surface $S$, the kernel is trivial, so $\Aut(S)$ embeds into the symmetric group of the ample cone.

\subsection{The roadmap}
We use a slightly different notation from \cite{Kondo1998}. For a subset $P$ in $NS(S)$, let $Sym(P)$ be the group of isometries of $NS(S)$ which keeps $P$ invariant. Element in $Sym(P)$ may not be realized as automorphisms, and we define $\Aut(P):=Sym(P)\cap \Aut(S)$.
Kond\=o's proof proceeds along the following steps:
\begin{enumerate}
	\item Each face $v$ of $D'$ is $r^\perp$ for some $r\in NS(S)$, the projection of some Leech root to $NS(S)$. For each $v$, identify an automorphism $f_v\in \Aut(S)$ which send one half-space divided by $v$ to the opposite half-space, or one of the half-spaces cut by another $r$. 
	\item Show that for each automorphism $f\in \Aut(S)$, there exists $h\in \{f_v\}_v$ such that  $h\circ f\in Sym(D')$. Then we conclude that $D'$ is a fundamental domain of $\Aut(S)$ on $\amp(S)$.
	\item Decide $\Aut(D')$. Then $\Aut(S)$ is generated by $\{f_v\}_v$ and $\Aut(D')$.
\end{enumerate}
We will prove Theorem \ref{Gstab} along the same steps, replacing $NS(S)$ with $A$, and $D'$ with a fundamental domain $\Omega$.

\subsection{Faces of the fundamental domain}
Recall $B:=\Z\{b, R\}$ and $A:=B^\perp$ with respect to $NS(S)$. Our goal is to prove 
\begin{align}
	\label{Omega}
	\Omega:=\overline{D'}\cap A
\end{align}
is a fundamental domain of $\Aut(S)_A$ acting on $\Nef(S)\cap A$.

Table \ref{fd} collects all the faces of $D'=D\cap P(S)$ and the classes $r$ such that $r^\perp$ cuts out the faces. We remark that classes here are scalar multiples of those from \cite{Kondo1998}, which does not change $r^\perp$. Those $N_\alpha$ and $T_\alpha$ in the table are $(-2)$-curves and they cut out the exterior wall of $D'$ along $\partial \Nef(S)$, so they do not correspond to any automorphisms. For each $r_w$, the associated $k_w$ appear as generators by \cite{Kondo1998}, and by \cite[Thm, 1.2]{Ohashi2009}, $k_w$ can be replaced by the HW-involutions $\varphi_w$.

\begin{table}[h]
	\centering
	\begin{tabular}[h]{r|c|c|l}\hline
		Count & $r$ & Type & Automorphism\\\hline
		$16$ & $N_\alpha$ & nodes & does not exist\\
		$16$ & $T_\alpha$ & tropes & does not exist\\
		$16$ & $\Lambda-2N_\alpha$ & projections & $p_\alpha$\\
		$16$ & $\sigma(\Lambda-2N_\alpha)$ & correlations& $q_\alpha$\\
		$60$ & $r_g:=\Lambda-\sum_{\alpha\in g} N_\alpha$ & HG-involutions & $\varphi_g$\\
		$192$& $r_w:=3\Lambda-2\sum_{\alpha\in w} N_\alpha$ & HW-involutions & $k_w$ or $\varphi_w$\\\hline
	\end{tabular}
	\caption{Faces of $D'$}
	\label{fd}
\end{table}

In particular:
\begin{align}
	\label{Dprime}
	D'=\{x\in \amp(x)\mid x\cdot r>0 \txt{for all $r$ in Table (\ref{fd})}\}.
\end{align}
We first note that $\Omega$ contains the projection $w''$ of $w'$ to $A\otimes \R$. Computation shows that 
\begin{align}
	\label{wpp}
w''=\frac{1}{7}\left(13\Lambda-3\sum_{\alpha\in J(C)_2} E_\alpha+4R\right).
\end{align}
In particular $\Omega$ is nonempty.
Recall that \[b={K_X}_{\mid S}\sim (3\Lambda-\sum_{i=1}^6 N_{i6}-2\sum_{1\leq i<j\leq 5} N_{ij})/2.\] (We use the convention that $N_{66}=N_0$.) Let
\begin{align}
	\label{cclass}
	c:=b+R=2\Lambda-\sum_\alpha N_\alpha.
\end{align}Then $c$ and $R$ also span $B$.
Since $R=T_0\in B$, every nef class $x\in F$ satisfies $x\cdot R=0$. This implies that $\Omega$ is in the face $R^\perp$ of $\Nef(S)$. Since $D'$ is rational polyhedral, so is $\Omega$. Therefore the faces of $\Omega$ are exactly those intersections of the faces $v$ of $D'$ with $A$, such that $v\cap A$ has maximal dimension.

We now exclude those faces of $v$ of $D'$ which do not cut a maximal dimensional face of $\Omega$. 
\begin{prop}
\label{dim10faces}
If $r$ in Table \ref{fd} associated to a projection, correlation, HG-involution of type 2, or HW-involution of type 2, then $\dim r^\perp\cap \Omega\leq 10$.
\end{prop}
\pf. Say $x$ is a $\Q$-nef class in $\Omega$. Then $x\cdot N_\alpha\geq 0$, so we can assume
\begin{align}
	x=a\Lambda-\sum_{\alpha\in J(C)_2} m_\alpha N_\alpha,
\end{align}
where $m_\alpha\geq 0$ and $a>0$. In particular we can assume $a=1$. Now $x\cdot r\geq 0$ for all $r$ in Table \ref{fd}, and $x\cdot R=x\cdot c=0$.   The equalities and inequalities read (coefficients reduced to integers):
\begin{alignat}{4}
	2&= \sum_{i=1}^6 m_{i6}, &&\quad (R)\label{eq1}\\
	4&= \sum_{\alpha} m_{\alpha}, &&\quad (c)\label{eq2}\\
	2&= \sum_{1\leq i<j\leq 5} m_{ij}, &&\quad (\ref{eq2})-(\ref{eq1})\label{eq3}\\
	2&\geq  \sum_{\beta\in I(T_\alpha)} m_{\beta}, &&\quad (T_\alpha)\label{eq4}\\
	1&\geq m_\alpha, &&\quad (\Lambda-2N_\alpha)\label{eq5}\\
	4&\geq  \sum_{\beta\not\in I(T_\alpha)} m_{\beta}, &&\quad (\sigma(\Lambda-2N_\alpha))\label{eq6}\\
	2&\geq  \sum_{\alpha\in g} m_{\alpha}, &&\quad (r_g)\label{eq7}\\
	3&\geq  \sum_{\alpha\in w} m_{\alpha}, &&\quad (r_w)\label{eq8}.
\end{alignat}
Now if $x$ is on the face cut by $r$, then the corresponding inequality will hold equal. We examine now each case. 

(i). When $r=(\Lambda-2N_\alpha)$. By symmetry we need only consider when $\alpha=0$ or $12$.

(ia). If $\alpha=0$, then $m_0=1$. We sum (\ref{eq4}) for $T_{i6}$, $i=1,\cdots,5$. Then
\begin{align}
10\geq 5+\sum_{i=1}^5 m_{i6}+2\sum_{1\leq i<j\leq 5} m_{ij}.
\end{align}
So $5\geq \sum_{i=1}^5 m_{i6}+2\sum_{1\leq i<j\leq 5} m_{ij}=5$ by (\ref{eq1}), (\ref{eq3}). Therefore the equality of (\ref{eq4}) holds for every $T_{i6}$, which together with (\ref{eq1}), (\ref{eq3}), gives $7$ linearly independent relations among $m_\alpha$. Hence the corresponding face has dimension at most $17-7=10$.

(ib). If $\alpha=12$, then $m_{12}=1$. Sum up  (\ref{eq4}) for $T_{\alpha}$, $\alpha=16,26,12,34,35,45$, and we have
\begin{align}
	12\geq 2\sum_{i=1}^6 m_{i6}+6+2\sum_{(ij)\neq 12} m_{ij}.
\end{align}
So $3\geq \sum_{i=1}^6 m_{i6}+\sum_{(ij)\neq 12} m_{ij}=2+1=3$. Therefore the equalities of (\ref{eq4}) hold for $\alpha=16,26,12,34,35,45$, and we have $7$ linearly independent relations.

(ii). When $r=\sigma(\Lambda-2N_\alpha)$, $4= \sum_{\beta\not\in I(T_\alpha)} m_{\beta}$. By (\ref{eq2}), we must have $m_\beta=0$ for all $\beta\in I(T_\alpha)$. This together with (\ref{eq1}) and (\ref{eq3}) gives at least $7$ linearly independent relations.

(iii). When $r=\Lambda-\sum_{\alpha\in g} N_\alpha$, we have $2=\sum_{\beta\in g} m_{\beta}$. If $g$ is of type 2, then by definition, $g\cap I(T_0)=\emptyset$. Now (\ref{eq3}) forces all $m_{ij}=0$ for $i,j\leq 5$, $ij\not\in g$. This together with (\ref{eq1}) and (\ref{eq3}), gives $8$ linearly independent relations.

(iv). When $r=3\Lambda-2\sum_{\alpha\in w} N_\alpha$, we have $3=\sum_{\beta\in w} m_{\beta}$. If $w$ is of type 2, then by definition, $\vv{w\cap I(T_0)}=1$. Suppose this common node is $i6$. Then
$m_{i6}\leq 1$ by (\ref{eq5}). Therefore 
\begin{align}
\sum_{\beta\in w,\, \beta\neq i6} m_{\beta}\geq 2.
\end{align}
Now (\ref{eq3}) forces all $m_{ij}=0$ for $i,j\leq 5$, $ij\not\in w$. This together with (\ref{eq1}) and (\ref{eq3}), gives $7$ linearly independent relations.\qed

\begin{prop}\label{15dim}
	The Cone $\Omega$ has maximal possible dimension $15$.
\end{prop}
\pf. We find a collection of linearly independent nef divisors in $\Omega$. On each face $r^\perp$ of $r=\Lambda-2N_{i6}$, there is a class
\begin{align}
C_{i}:=5\Lambda-5N_{i6}-\sum_{\alpha\neq i6} N_\alpha.
\end{align}
For each G\"opel tetrad $g$ of type $2$, the class 
\begin{align}
	F_{g}:=6\Lambda-2\sum_{i=1}^6 N_{i6}-3\sum_{\alpha\in g} N_\alpha
\end{align}
is on the face $r_g^\perp$.
It is clear that $C_i$ and $F_g$ all satisfy all the inequalities and equalities (\ref{eq1}) to (\ref{eq8}). Hence they are in $\Omega$.
Then direct calculation shows that the matrix with rows $\{C_i, F_g\}$ has rank $15$. Hence $\dim \Omega=15$.\qed

As a conclusion, the faces of $\Omega$ consists of at most those $r^\perp$ where $r$ is one of $T_\alpha$, $N_\alpha$, $r_g$ of type $1$ and $r_w$ of type $1$.

\subsection{The generation of $\Aut(S)_A$}
Define $G$ as the group generated by all those HG-type automorphisms $z_g$ and Keum's automorphisms $z_w$.

\begin{lemma}
	\label{subgp}
	The group $G$ is a subgroup of $\Aut(S)_A$.
\end{lemma}
\pf. We need only show all $z_g$ and $z_w$ fixes $A$. The definition of $z_g$ and $z_w$ implies that they fix $R$. They fixes $c$ since among the generators of $\Aut(S)$, the translations and HG-involutions fix $c$, while the switch and the HW-involution send $c$ to $-c$.\qed

\begin{lemma}
	\label{pairingfaceauto}
	For each $z_g$,  $z_g(r_g)=-r_g$ and $z_g$ exchanges the two half-spaces divided by $r_g$.  Each $z_w$ exchanges the two half-spaces $w_{>0}$ and $w'_{<0}$, where $w'$ is the dual of $w$.
\end{lemma}
\pf. $z_g(r_g)=-r_g$ is a direct calculation. An alternative proof is like this: recall that $z_g=\varphi_g\circ t_\alpha$ for some $\alpha$. In either case $t_\alpha(r_g)=r_g$, and $\varphi_g$ acts on $r_g$ as a reflection about the itself, so $\varphi(r_g)=-r_g$. Consider a class $x=a\Lambda-\sum_{\alpha\in J(C)_2} m_\alpha N_\alpha$. The same calculation shows that 
\begin{align}
	z_g(x)\cdot r_g=-4a+2\sum_{\beta\in g} m_\beta=-x\cdot r_g,
\end{align}
so $z_g$ exchanges the two half-spaces.

For $w$ of type 1, by definition, $z_w=\varphi_w\circ \sigma\circ t_\alpha$. In either case, it can be computed that $(\sigma( t_\alpha(x))\cdot r_w=x\cdot r_{w'}=12a-4\sum_{\beta\in w'} m_\beta$. On the other hand, $\varphi_w(x)\cdot r_{w}=-12a+4\sum_{\beta\in w} m_\beta=-x\cdot r_{w}$. Therefore
\begin{align}
	z_w(x)\cdot r_{w}=-x\cdot r_{w'}.
\end{align}
\qed

Now we have proved (1) of the roadmap. We need to prove (2).
\begin{prop}\label{homing}
	For each automorphism $f\in \Aut(S)_A$, there exists $h\in G$ such that  $h\circ f\in Sym(\Omega)$.
\end{prop}
\pf. The method is from \cite[7.3]{Kondo1998}. The projection $w''$ (\ref{wpp}) of $w'$ to $A$ lies in $\Int(\Omega)$, the interior of $\Omega$. Let $u=f(w'')$. Let $h\in \{z_g, z_w\}$ be the automorphism which achieves the minimum of $\{h(u)\cdot w''\mid h\in G\}$. By Lemma \ref{subgp}, $G$ fixes $A$, so $h(u)=h(f(w''))\in A$.
Therefore we show $h(u)\in \Int(\Omega)$. The non-nef-boundary faces of $\Omega$ correspond to $w$ or $g$ of type 1. For any face $v$ of $\Omega$, let the corresponding automorphism be $y_v$. Then
\begin{align}
	h(u)\cdot w''\leq y_v\ik(h(u))\cdot w''=h(u)\cdot y_v(w'').
\end{align}
There are two cases.

(i). If $y_v=z_g$, then it can be computed
\begin{align}
	z_g(w'')=w''+2r_g.
\end{align}
Therefore 
\begin{align}
	h(u)\cdot w''\leq h(u)\cdot y_v(w'')=h(u)\cdot w''+2h(u)\cdot r_g.
\end{align}
Therefore $h(u)\cdot r_g\geq 0$.

(ii). If $y_v=z_w$, then \begin{align}
	z_w(w'')=w''+3r_w.
\end{align}
Therefore $h(u)\cdot r_w\geq 0$. This proves that $h(u)\in \Omega$. Now for any $f'\in \Aut(S)_A$, either $f'$ fixes $\Omega$ or $f'(\Int(\Omega))\cap \Int(\Omega)=\emptyset$. Since $w\in \Int(\Omega)$, we must have that $h(u)\in \Int(\Omega)$ too. Therefore $h\circ f\in \Aut(S)_A$ must fix $\Omega$, which finishes the proof.\qed

\begin{lemma}
	\label{unique-2}
	The only $(-2)$-curve class in $B$ is $R=T_0$.
\end{lemma}
\pf. Suppose $x=s c+t R$ is a $(-2)$-curve class. Recall (\ref{cclass}). In order for $s c+tR\in NS(S)$, we must have $s\in (1/2)\Z$ and $t\in \Z$. Then we can assume $s=v/2$ and $v\in \Z$. Then
\begin{align}
	-2=\left(\frac{vc}{2}+tR\right)^2=-2(2v^2+tv+t^2).
\end{align}
The only integer solutions are $v=0,t=\pm 1$. That is, $x=\pm T_0$.\qed
\begin{lemma}
	\label{symdream}
	The group $\Aut(\Omega)=\{\id\}$.
\end{lemma}
\pf. The group $Sym(\Omega)$ is a subgroup of $\Aut(D')$. Indeed, by \cite[\S 4]{Kondo1998}, every isometry $f$ of $NS(S)$ extends to an isometry of $\II_{1,25}$, which either fixes $D'$ or sends $D'$ to another fundamental domain $D''$. Observe that $A\subset R^\perp$ is not a face of $R^\perp\cap \Nef(S)$, and $\Omega = \overline{D'}\cap A$ is not a face of $D'$, because neither $\pm b$ is effective, or scalar multiple of the those $r$ that cut out the faces of $D'$ (Table (\ref{fd})). 

Suppose now $f$ fixes $\Omega$ and $f(D') = D''$ is not $D'$. Then $\overline{D'} \cap \overline{D''}$ must be a facet of $D'$. However $\Omega\subset \overline{D'} \cap \overline{D''}$. Therefore $\overline{D'} \cap \overline{D''}$ must be the face $D'\cap R^\perp$ since $\dim \Omega = 15$. On the other hand, $R$ is effective, so $R^\perp$ is a facet of $\Nef(S)$, so both $D'$ and $D''$ lie in the half plane $R_{>0}$. Hence $D'\cap D'' \cap R^\perp$ can't be a facet (namely, a codimension $1$ face) of $D'$, which is a contradiction. Therefore $f$ must fix $D'$.

By \cite[7.2]{Kondo1998}, $\Aut(D')$ is a $2$-elementary abelian group isomorphic to $(\Z/2\Z)^5$, generated by the $16$ translations $t_\alpha$ and the switch $\sigma$ (the permutations of the six Weierstrass points do not realize as automorphisms of $S$). A simple calculation shows that no nontrivial automorphisms in $\Aut(D')$ fixes $T_0$. By Lemma \ref{unique-2}, we conclude that they do not fix $B$, and hence $A$. Therefore $\Aut(\Omega)$ is trivial.\qed

\pfof{Theorem \ref{Gstab}}. We proved (1) (2) of the roadmap above, so $G$ acts on $\Nef(S)\cap A$ with a fundamental domain $\Omega$, where we can find $45+120$ faces corresponding to $z_g$ or $z_w$, and the rest are cut out by some of the $N_\alpha$ and $T_\alpha$.
Now Proposition \ref{homing} and Lemma \ref{symdream} implies that $G=\Aut(S)_A$.

Finally by definition $\Aut(S)_b\cap \Aut(S)_R\subset \Aut(S)_A$. Conversely it follows from that every generator of $G$, namely, those $z_g$ and $z_w$ where $g$ and $w$ are of type 1, fixes $R$ and $c$, hence also $b = c -R$. \qed

\bibliography{mybib.bib}
\bibliographystyle{alpha}
\end{document}